\documentclass[11pt]{article}
\usepackage[english]{babel}
\usepackage{amsmath,latexsym}
\usepackage[psamsfonts]{amssymb}
\usepackage[dvips]{graphicx}
\usepackage[mathcal]{euscript}
\begin{document}
 515.168
\begin{center}
{\bf On the geometry of the orbits of Killing vector fields } \\[2mm]
{\bf  A.Ya. Narmanov J. O. Aslonov}\\
\end{center}\vspace{5mm}

\begin{abstract} Let $D$  be a set of smooth vector fields on the smooth
manifold $M$.It is known that orbits of $D$ are submanifolds of M.
Partition $F$ of M into orbits of $D$ is a singular foliation. In
this paper we are studying geometry of foliation which is
generated by orbits of a family of Killing vector fields.In the
case $M=R^3$ it is obtained full geometrical classification of
$F$. Throughout this paper the word "smooth" refers to a class
$C^\infty$.
\end{abstract}

Let $M$ is a smooth, connected Riemannian manifold of dimension $n$,
$X$ is a smooth vector field, $ X ^ t (x) $ is an integral curve passing
through $ x $ for $ t = 0 $.

\textbf{Definition 1.} Vector field $ X $ on $ M $ is called a
Killing field if the infinitesimal transformation $ x \to X ^ t
(x) $ is an isometry of $M$ for any $t$.

\textbf{Example 1.} In the three-dimensional Euclidean space $R^3
(x, y, z) $, there are six linearly independent Killing fields on
the field of real numbers:

$$ X_1 = \frac{\partial}{\partial x}, \ X_2 = \frac {\partial}
{\partial y}, \ X_3 = \frac {\partial } {\partial z},
X_4 =-y \frac {\partial} {\partial z} + z \frac {\partial}
{\partial y}, \ X_5 =-z \frac {\partial} {\partial x} + x \frac
{\partial} {\partial z}, \ X_6 =-x \frac {\partial }
{\partial y} + y \frac {\partial} {\partial x}. $$

The groups of transformations generated by vector fields $ X_1, \
X_2, \ X_3 $ are groups of translations in the direction of the
axes $ Ox, \ Oy $ and $ Oz $, respectively, and the groups of
transformations generated by   last three vector fields are
rotations around the axes of the $ Ox, \ Oy $ and $ Oz $
accordingly.

The last three fields are Killing fields on the sphere $ S ^ 2 $
too.

\textbf{Example 2.} Consider the three-dimensional sphere $ S ^ 3
$ in $ R ^ 4 \thickapprox C ^ 2 $ with the induced metric. Let $
(x_1, \ x_2, \ x_3, \ x_4) $ is a point on the sphere $ S ^ 3 $.
With complex numbers $ S ^ 3 $ can be described as follows: $ S ^
3 = \{(z_1, \ z_2); \ | z_1 | ^ 2 + | z_2 | ^ 2 = 1 \} $ where, $
z_1 = x_1 + ix_2 $, $ z_2 = x_3 + ix_4 $.

Consider in $ R ^ 4 $ Killing vector field
$$ X =-x_2 \frac {\partial}
{\partial x_1} + x_1 \frac {\partial} {\partial x_2}
- x_4 \frac {\partial} {\partial x_3} + x_3 \frac {\partial}
{\partial x_4}. $$

It is easy to check that this vector field is tangent to the
sphere. For a point $ (z_1, \ z_2) \ in S ^ 3 $ the integral curve
of the vector field $ X $, starting from the point $ (z_1, \ z_2)
$ for $ t = 0 $ has the form

$$ \gamma (t) = \{(z_1e ^ {it}, z_2e ^ {it}), \ - \infty <t <\infty \}. $$

It is obvious that the integral curve $ \gamma (t) $ is circle.
The family of integral curves of the vector field $ X $ generates
a smooth bundle, which is called Hopf bundle .

In the future, will require the following statement \cite{Berest}

\textbf{Proposition.} Vector field
\[X = \sum_ {i = 1} ^ n \xi _i
\frac {\partial} {\partial x_i} \]
in $ R ^ n $ is a Killing field
if and only if the conditions are satisfied

\[\frac {\partial \ xi _i} {\partial x _j} + \frac {\partial
\xi _j} {\partial x _i} = 0 \ i \neq j, \ \frac{\partial \xi _i}
{\partial x _i} = 0, \ i = 1, ..., n. \]

It is known that if the length of the Killing vector field is
constant on the whole manifold, the integral curves are geodesics
lines \cite{arnold}.

The following lemma show that on  two-dimensional circular
cylinder, the integral curves of Killing vector field is always
are geodesics.

\textbf{Lemma 1}. Each integral curve of a smooth Killing vector
fields on the two-dimensional circular cylinder are geodesic
lines.

\textbf{Proof.} Let the two-dimensional circular cylinder $ M $
parameterized as follows
\[\left \{\begin{array} {ll}
x = sin u\\
y = cos u \\ z = v \\
\end{array} \right. \]

In three-dimensional Euclidean space $ R ^ 3 $, consider the
following Killing vector fields

$$ X_1 =-x \frac {\partial} {\partial y} + y \frac {\partial} {\partial x},
 X_2 = \frac {\partial} {\partial z}. $$

It is easy to verify that these vector fields tangent to $M$, and
the are linearly independent vector Killing fields on $M$.

It is easy to verify that the integral curves of these vector
fields are are geodesic lines on $M$.

Let $ X $ is a smooth Killing vector field on $ M. $ Then the there
are smooth functions of $ \lambda _1 (x, y, z) $ and $ \lambda _2
(X, y, z)$ such that the vector field $ X $ has the form $$ X =
\lambda _1 (x, y, z) X_1 + \lambda _2 (x, y, z) X_2. $$

Since the vector field $ X $ is the Killing vector field, applying
proposition to the vector field $ X, $ we have the following:
\begin{equation} \lambda _1 (x, y, z) = \lambda _1
(z), \quad \lambda _2 (x, y, z) = \lambda _2 (x, y), \quad \frac
{\partial \lambda _1} {\partial z} + y \frac {\partial \lambda _2}
{\partial x} = 0, \quad \frac {\partial \lambda _1} {\partial z}-x
\frac {\partial \lambda _2} {\partial y} = 0.  \end{equation}

It is clear that the function $ \lambda _1 (x, y, z) $ does not
depend on variables $ x, y, $ and the function $ \lambda _2 (x,
y, z) $ does not depend on variable $ z. $

It is known that the Lie bracket of two Killing vector fields is
Killing field [2]. Thus, the vector fields $ [X, X_1] $ and $ [X,
X_2] $ are  Killing vector fields. The vector field $ [X, X_2] $
is given by $ [X, X_2] = \lambda _1 [X, X_2] + X_2 (\lambda _1)
X_1 + \lambda _2 [X_2, X_2] + X_1 (\lambda _2) X_2. $

Since $ [X_1, X_2] = 0 $ we have $ [X, X_2] = X_2 (\lambda _1) X_1
+ X_1 (\lambda _2) X_2. $

Now, applying proposition to the vector field $[X, X_2],$ we
obtain that $ \frac {\partial \lambda _1} {\partial z} = 0. $
Hence, using the that $ (1) $ we obtain $ \frac {\partial \lambda _2}
{\partial x} = \frac {\partial \lambda _2} {\partial y} = 0. $

Thus, $ \lambda _1 (x, y, z) $ and $ \lambda _2 (x, y, z) $ are
constant functions.

To find the integral curve of the vector field $X$ consider the
system of differential equations: \[\left \{\begin{array} {ll}
\frac {\partial x} {\partial t} = \lambda _1 y \\ \frac{\partial
y} {\partial t} = - \lambda _1x \\ \frac {\partial z}
{\partial t} = \lambda _2 \\
\end{array} \right. \]

The solution of this system with initial condition $ x (0) = x_0,
y (0) = y_0, z (0) = z_0 $ is given by

\begin{equation} \left \{\begin {array} {ll} x (t) = x_0 cos
\lambda _1 t + y_0 sin \lambda _1 t \\ y (t) =-x_0 sin \lambda
_1 t + y_0 cos \lambda _1 t \\ z (t) = \lambda _2 t + z_0. \\
\end{array} \right. \end{equation}

If the point $ (x_0, y_0, z_0) $ belongs to the cylinder, then $
x_0 ^ 2 + y_0 ^ 2 = 1 $ and equation (2) defines a geodesic line
on the cylinder, which is a helical line, if $\lambda _1 \neq 0,
\lambda _2 \neq 0 $. If $ \lambda _1 = 0, \lambda _2 \neq 0 $,
then the geodesic is a straight line, if  $ \lambda _1 \neq 0,
\lambda _2 = 0 $, then equation (2) defines a circle.

\textbf{Remark.} The next example show, the integral curves of the
Killing vector field in three-dimensional circular cylinder must
not be geodesics.

\textbf{Example 3.} Let, $ M = S ^ 2 \times R ^ 1 = \{(x, y,
z, w) \in R ^ 4 : X ^ 2 + y ^ 2 + z ^ 2 = 1 \}.$

Consider the vector field $X = y \frac {\partial} {\partial x}-x
\frac {\partial} {\partial y} $ in $ R ^ 4. $ By proposition, we
can check that this field is the Killing field in $ R ^ 4. $
Furthermore, this vector field is tangent to $ M $ and it is
Killing vector field on $ M. $ The integral curve of $ X $ passing
through the point $ (x_0, y_0, z_0, w_0) $  has the
form \[\left \{\begin{array} {ll} x (t) = x_0 cos t + y_0 sin t \\
y (t) = - x_0 sin t + y_0 cos t \\ z (t) = z_0 \\ w (t) = w_0 \\
\end{array} \right. \]

If $ z_0 \neq 0 $ then this curve is not a great circle on $ S ^
2. $ Hence, it is not a geodesic.

\textbf{Theorem 1. }Let $D$ is a family of smooth Killing vector
fields
 on $M.$ We assume that the dimension of the every
orbit of family $D$ smaller than $n.$ Then the partition of the
manifold $M $ into the orbits  is a singular Riemannian foliation.

\textbf{Proof.}Let us to recall that   foliation is riemannian if
every geodesic that is perpendicular at one point to a leaf
remains perpendicular to every leaf it meets.

Let $F$ is a partitionn of the manifold into the orbits of the
family $D.$ As a result of the work of Stefan and  Sussmann $ F $
is a singular foliation. As follows from the results of Molino
\cite{molino} and A. Narmanov \cite{narmanov1996}, the foliation $
F $ is Riemannian.

Let $D$ be a set family of smooth vector fields on manifold
 $M$.The family $D$ generates a distribution $P_D:x\to
 P_D(x)$,where $x\in M$,$P_D(x)$is the linear hull of the vector$X(x)$,$X\in D.$
We call family $D$ completely integrable if distribution $P_D$ is
completely integrable.

Consider the vector fields $ X $ and $ Y $ in $ R ^ 4 $ are in
Cartesian coordinates are as follows: $$ X = x \frac {\partial}
{\partial y}-y \frac {\partial} {\partial x} +
z \frac{\partial} {\partial w}-w \frac {\partial} {\partial z}, \\ Y =
z \frac {\partial} {\partial x}-x \frac {\partial} {\partial
z} + w \frac {\partial} {\partial y}-y \frac {\partial} {\partial w}. $$

By proposition, it is easy to verify that these fields are Killing
fields. Consider the unit sphere $ S ^ 3 $ induced metric in $ R ^
4 $ is given by the equation $ x ^ 2 + y ^ 2 + z ^ 2 + w ^ 2 = 1.
$

It is easy to check that these vector fields are tangent to the
sphere $ S ^ 3. $ Hence, these fields are Killing vector fields on
a sphere $ S ^ 3.$

\textbf{Theorem 2.} The family of vector fields $ X, Y $ is
completely integrable.The decomposition of $ S ^ 3$ to the orbits
is a singular Riemannian foliation, regular leaves of
 which are two-dimensional tories. The set of singular leaves consists of two circles.

\textbf{Proof.} With easily calculations can be  showen that $[X,
Y] = 0.$ Hence, that the family $D =\{X, Y \}$ is involutive, and
by   Frobenius-Hermann theorem, family of vector fields $D = \{X,
Y \}$ is completely integrable.

We now show that each regular leaf of foliation $F$ generated by
orbits of family $ X, Y $ are two-dimensional torus.

Consider the corresponding system of differential equations
\begin{equation} \left \{\begin{array} {ll} \frac {\partial
x} {\partial t} =-y \\
\frac {\partial y} {\partial t} = x \\
\frac {\partial z} {\partial t} = - w \\
\frac {\partial  w} {\partial t} = z \\
\end{array} \right. \qquad \qquad \qquad
\left \{\begin{array} {ll}
\frac {\partial x} {\partial t} = z \\
\frac {\partial y} {\partial t} = w \\
\frac {\partial z} {\partial t} =-x \\
\frac {\partial w} {\partial t} =-y. \\
\end{array} \right.
\end {equation}

For a point on the sphere integral curve $ \gamma _ {0} ^ {1} $
vector fields $ X, $ passing through $ p_0 $ for $ t = 0 $, has
the following parametric equations

\[\left \{\begin{array} {ll}
x (t) = x _0 cos t - y_0 sin t \\
y (t) = x_0 sin t + y_0 cos t \\
z (t) = z_0 cos t - w_0 sin t\\
w (t) = z_0 sin t + w_0 cos t \\
\end{array} \right. \]

It is easy to see that this is a closed curve.

The integral curve $ \gamma _ {0} ^ {2} $ vector field $ Y, $
passing through the point $ p_0 $ for $ t = 0 $, given by the
following parametric equations

\[\left \{\begin{array} {ll}
x (t) = x _0 cos t + z_0 sin t \\
y (t) =-y_0 cos t + w_0 sin t \\
z (t) = x _0 sin t + z_0 cos t \\
w (t) = y _0 sin t + w_0 cos t \\
\end{array} \right. \]

The integral curve $ \gamma _ {0} ^ {2} $ is also closed.

The vector fields $ X $ and $ Y $ are collinear only at points of
two circles, which are the intersections of spheres with
two-dimensional planes
\begin{equation} \left \{\begin{array} {l}
x = w, \\ y =-z \\ \end{array} \right. \qquad \qquad \left \{\begin{array}{ll}
x =-w, \\ y = z \\ \end{array} \right. \end{equation}
respectively. These circles are given by equations
\begin{equation} \left \{\begin{array}{ll}
x = w, \\ y =-z \\ 2x ^ 2 +2 y ^ 2 = 1 \\
\end{array} \right. \qquad \left \{\begin{array} {ll}
x =-w, \\ y = z \\
2z ^ 2 +2 w ^ 2 = 1. \\ \end{array} \right.
\end{equation}

We show that these circles are integral curves for vector fields $
X $ and $ Y. $

Let a point $ p_0 (x_0, y_0, z_0, w_0) $ belongs to the first
circle. Then $ x_0 = w_0, \ y_0 =-z_0 $ and integral curve of the
vector field $ X, $ issuing from $ p_0, $ has the form:

\[\left \{\begin{array} {ll} x (t) = x _0 cos t - y_0 sin t \\
y (t) = x_0 sin t + y_0 cos t \\ z (t) =-y _0 cos t - x_0 sin t\\
 w (t) =-y _0 sin t + x_0 cos t. \\ \end{array} \right. \]

It is easy to see that $ x (t) = w (t), \ y (t) =-z (t), $ for all
$ t \in R. $ This means that a circle is an integral curve for
vector field $ X $.

Now consider the integral curve of the vector field $ Y, $
extending from a point $ p_0. $ It is given by the following
parametric equations: \[\left \{\begin{array} {ll}
x (t) = x_0 cos t - y_0 sin t \\ y (t) =-y _0 cos t + w_0 sin t \\
z (t)=-x _0 sin t - y_0 cos t \\ w (t) =-y _0 sin t + x_0 cos t. \\
\end{array} \right. \]

It is easy to see the $ x (t) = w (t), \ y (t) =-z (t), $ for all
$ t \in R. $ This means that a circle is an integral curve for
vector field $ Y $. Similarly we can prove that the second circle
is also an integral curve for the vector fields $ X $ and $ Y. $

Now consider a point $ p $ on the sphere, which does not belong to
the above considered e circles. Denote by $ \gamma _1 (t) $
integral curve of the vector field $ X, $ issuing from the point $
p $. At points $ \gamma _1 (t) $ vectors $ X $ and $ Y $ are
linearly independent.

Consider the mapping $ (t, s) \to X ^ t (Y ^ s (p)), $ for all $
(t, s) \in R ^ 2. $

It is clear that the image of two-dimensional plane $ R ^ 2 $ is
the orbit family $ D $ passing through $p. $ Since the vectors $ X
$ and $ Y $ are linearly independent at a point $ p, $ rank  of
this map is equal to two. Thus, the orbit passing through $ p $ is
two-dimensional manifold.

The equality $ [X, Y] = 0 $ means that $ X ^ t (Y ^ s (p)) = Y ^ s
(X ^ t (p)), $ for all $ (t, s) \ in R ^ 2. $ This means that the
flow of the vector field $ X $ maps the integral curves of the
vector field $ Y $ to the integral curves of the vector field $ Y
$ (respectively, the flow of the vector field $ Y $ maps the
integral curves of vector fields in integral curves of the vector
field $ X $).

Since the integral curve $ \gamma _2 (s) = Y ^ s (X ^ t (p)) $ of
vector field $ Y, $ issuing from $ \gamma _1 (t) = X ^ t (p) $ at
$ s = 0 $ is a closed curve, if we consider all integral curves $
\gamma _2 (s) = Y ^ s (X ^ t (p)) $ for every point of $ \gamma _1
(t) = X ^ t (p), $ we obtain the two-dimensional torus.

Both vector fields $ X $ and $ Y $ tangent to the torus. As
follows from Lemma \cite{olver} if vector fields $ X $ and $ Y $
is everywhere tangent to the submanifold $ N $, then their Lie
bracket $ [X, Y]$ is tangent to the manifold $ N $.

Therefore, while we move on this torus, we can not leave it,
therefore, the orbit of our family passing through $ P $ is
contained on the torus. So as, from a point $ P $, we can come to
any point of the torus, moving along the integral curves of vector
fields $ X $ and $ Y, $ the orbit of the family $ D $ containing
the point $ P, $ is two-dimensional torus. The fact that the
orbits generate Riemannian foliation follows from Theorem
1.Theorem 2 is proved.

\textbf{Remark.} Vector field $ X $ and $ Y, $ does not have
critical points. Since each vector field tangent to the
two-dimensional  sphere must have a critical point on it, the
orbit of the family can not be two-dimensional sphere
\cite{arnold}.

The following theorem gives a complete classification of foliation
$F$ generated by orbits of family of Killing vector fields in
three-dimensional Euclidean space.

\textbf{Theorem 3.} Let $ D $ - a family of Killing vector fields
in $ R ^ 3. $ Then the orbits of this family generate a foliation
$ F, $ which is one of the following seven types:

1) foliation $ F $ consists of parallel lines;

2) foliation $ F $ consists of concentric circles, lying on the
parallel planes and straight line,which is a set of centers.

3) foliation $ F $ consists of a helical lines lying on
concentric circular cylinders;

4) foliation $ F $ consists of parallel planes;

5) foliation $ F $ consists of concentric spheres and a point
(the center spheres);

6) foliation $ F $ consists of concentric circular cylinders
and line (the axis of the cylinder);

7) foliation $ F $ has only one leaf $ R ^ 3. $

 We will use the following lemma which was proved in the paper
 \cite{narmanov1996}.

\textbf{ Lemma.2.} Let $ F $ - singular Riemannian foliation on
complete Riemannian manifold $ M $, $ \gamma _0 $ - geodesic going
from some point $ x_0 $ to some point $ y_0 $ orthogonal to $ F $.
Then for each point $ x \in L (x_0) $ there exists a geodesic
  $\gamma $, going from $ x $ to some point of
the leaf $ L (y_0)$,orthogonal to $ F $ and length of $\gamma$
equal to length of $ \gamma _0 $.

\textbf{Proof.}

1) Consider the case when there exists a unique point $ p_0 $ such
that $ L (p_0) = {p_0} $. This means that the isometries generated
by vector fields in $ D $, are the rotations around the axes
passing through the point $ p_0 $. In addition, the number of
vector fields is greater than one.

Let $ S_r ^ 2 $ - sphere of radius $ r> 0 $ with center at point $
p_0 $. Then isometries generated by vector fields from $ D $, map
 $ S_r ^ 2 $ to $ S_r ^ 2 $ .
In addition, since we have more than one rotation, orbit $ L (q)
$, passes through the point $ q \in S_r ^ 2 $ is two-dimensional
manifold. Since we have a rotation around non-parallel axes
passing through the point $ p_0 $, the orbit
 $L (q) $ coincides with the manifold $ S_r ^ 2 $. In this case, the
foliation $ F $ consists of concentric spheres and a point $p_0$
(the center of the spheres);

2) There are two different points $ p_1 $ and $ p_2 $ such that $
L (p_1) = {p_1} $, $ L (p_2) = {p_2} $.

This means that an isometries generated by vector fields from  $ D
$, remain fixed points $ p_1 $ and $ p_2 $. It follows that all
points of straight line $ p_1p_2 $ are fixed. In this case, we
have no parallel transports, we have only rotation around a
straight line $ p_1p_2 $. Hence, for all points $qQ $, do not
belong the straight $ p_1p_2 $, the orbit is a circle with
centered on the line $ p_1p_2 $. Foliation $ F $ consists of
concentric circles lying in parallel planes perpendicular to line
$ p_1p_2 $ and line  $ p_1p_2 $.

3) Let $ dimL (p) = 1 $ for all $ p \in R ^ 3 $. In this case
there is no fixed points. Therefore, the vector fields from $ D $
does not generate only the rotation, they may consist of parallel
transports and the composition of parallel transports and
rotations. If we have only parallel transports, all the orbits are
parallel lines. Foliation $ F $ consists of parallel lines.

If there is a composition of a translation and rotation, then by
Lemma-1 each orbit, except one, are helical lines lying on
concentric cylinders with a common axis. Foliation $ F $ consists
of helical lines lying on concentric circular cylinders, one of
which is the axis of the cylinder.

4) Let $ dimL (p) \geq 1 $, and there is a point $p$, for which $
dimL (p) = 2 $ and $ dimL (p) = 1 $.

Let $ dimL (p_0) = 1 $ for the point $ p_0 \in R ^ 3 $. There are
two features:

a) An orbit $ L (p_0) $ - straight line

or b) an orbit $L (p_0) $ - helix.

Consider the case when $ L (p_0) $ -straight line. In this case,
the isometries generated by vector fields from $ D $, can not
consist of translations only . There is a composition of
translations and rotations about the line $ L (p_0) $.
Consequently, we have a family of concentric cylinders with common
axis $ L (p_0) $. Foliation $ F $ consists of concentric circular
cylinder and a straight line (axis cylinders).

Consider the case where $ L (p_0) $ - helix.In this case there is
a point $ p_1 $, such that the orbit  $ L (p_1) $ is a straight
line. Therefore, this case reduces to the previous one.

5) All the orbits are two-dimensional manifolds.In this case we
have a regular two-dimensional Riemannian foliation. The results
of \cite{narmanov2011} implies that the foliation consists of
parallel planes.

6) There is a three-dimensional orbit of $ L (p_0) $. In this
case,by the Theorem of H. Sussmann \cite{sussman} and P.Stefan
\cite{stefan}, the orbit of $ L (p_0) $ is manifold, and therefore
all its points are interior points.

Let $ \varepsilon> 0 $ such that the open ball $ B_ {\varepsilon}
(P_0) $ of radius $ \varepsilon $ centered at $ p_0 $ is contained
in $ L (p_0) $. We show that $ L (p_0) $ is a closed set in $ R ^
3 $. Let $ p_m \in L (p_0) $, and $ p_m \to q $ for $ m \to \infty
$. Then, for sufficiently large $ m $ we have $d (p_m, q)
<\varepsilon $, where $ d (p_m, q) $ is  the distance in $ R ^ 3 $
between points $ p_m $ and $ q $. Due to the fact that $ p_m \in L
(p_0) $ exists isometry $ f_m $, generated by vector fields from $
D $ such that $ f_m (p_0) = p_m $. Then $ f_m ^ {-1} (p_m) = p_0
$, $ f_m ^ {-1} (q) \in B_ {\varepsilon} (p_0) $. This implies
that $ q \in L (p_0) $, i.e. $ L (p_0) $ is a closed set. By the
connectedness of $ R ^ 3 $ orbit of $ L (p_0) $ coincides with $ R
^ 3 $.

\end{document}